\documentclass[11pt, twoside, reqno]{article}
\usepackage{amsmath,amsthm}
\usepackage{amssymb,latexsym}
\usepackage{enumerate}
\usepackage{amsfonts}
\usepackage{mathrsfs}
\usepackage{amsmath, amsthm}
\usepackage{amssymb}
\usepackage{fancyhdr}

\makeatletter

\@addtoreset{equation}{section}

\makeatother

\newfam\msbfam
\font\tenmsb=msbm10    \textfont\msbfam=\tenmsb
\font\sevenmsb=msbm7 \scriptfont\msbfam=\sevenmsb
\font\fivemsb=msbm5 \scriptscriptfont\msbfam=\fivemsb

\newfam\bigfam
\font\tenbig=msbm10 scaled \magstep2   \textfont\bigfam=\tenbig
\font\sevenbig=msbm7 scaled \magstep2 \scriptfont\bigfam=\sevenbig
\font\fivebig=msbm5 scaled \magstep2
\scriptscriptfont\bigfam=\fivebig

\newtheorem{thm}{Theorem}[section]
\newtheorem{lem}{Lemma}[section]

\begin{document}
	

\renewcommand{\thefootnote}{\fnsymbol {footnote}}
	\title{{\bf Weighted boundedness of   Hardy and P\'{o}lya-Knopp type operators with variable limits on product space}}
	\author{{\bf Qianjun HE\footnote{{Corresponding author}}  \quad \bf Dunyan YAN}
		\\\footnotesize\textit{School of Mathematical Sciences, University of Chinese Academy of Sciences, Beijing 100049, P. R. China}}
\footnotetext {{Supported by the National Natural Science Foundation of China (grant numbers 11471309 and 11561062).}}
\footnotetext {{E-mail address: heqianjun16@mails.ucas.ac.cn (Qianjun HE); ydunyan@ucas.ac.cn (Dunyan YAN)}}
		\date{}
		\maketitle

		\noindent{{\bf Abstract}} \quad We  characterize the sufficient conditions which three weight functions $u$ and $v_{1}, v_{2}$  satisfy  ensure the boundedness of  the Hardy operator with variable limits on product space. The corresponding bound  is explicitly worked out. Moreover, as application,  we can obtain  an explicit  scale of bound for the P\'{o}lya-Knopp type operator with  certain weights.
		\medskip
		
		\noindent{{\bf MR(2000) Subject Classification}} \quad 42B20, 42B25

		
		{\centering\section{ Introduction }}

	\thispagestyle{empty}
	We study the Hardy type operator defined as
		\begin{equation}\label{eq:1}
		(H_{2}f)(x_{1},x_{2})=\int_{a_{1}(x_{1})}^{b_{1}(x_{1})}\int_{a_{2}(x_{2})}^{b_{2}(x_{2})}f(t_{1},t_{2})dt_{1}dt_{2},
		\end{equation}
		where both functions $a_{i}$ and $b_{i}$ are strictly increasing differentiable functions on $\mathbb{R}_{+}=(0,\infty)$, and  satisfy
		\begin{equation}\label{eq:2}
		{
			\begin{split}
			&a_{i}(0)=b_{i}(0)=0,\\
			&a_{i}(x_{i})<b_{i}(x_{i}),\\
			&a_{i}(\infty)=b_{i}(\infty)=\infty,
			\end{split}
			}
		\end{equation}
for  $0<x_{i}<\infty$ with $i=1,2$.

		For one dimensional case,  for nonnegtive measurable $f\in L^p(\mathbb{R}_{+})$ with $1<p\leq q<\infty$,  Ushakova \cite{U}  obtained that
		\begin{equation}\label{eq:3}
		\left(\int_{0}^{\infty}\left(\int_{a(x)}^{b(x)}f(t)dt\right)^{q}u(x)dx\right)^{\frac{1}{q}}\leq C\left(\int_{0}^{\infty}f^{p}(x)v(x)dx\right)^{\frac{1}{p}}
		\end{equation}
		holds if and only if
 		\begin{equation}\label{eq:4}
		B(s)=\sup\left(\int_{t}^{x}u(y)\left(\int_{a(y)}^{b(y)}v^{1-p^{\prime}}(\eta)d\eta\right)^{\frac{q(p-s)}{p}}dy\right)^{\frac{1}{q}}\left(\int_{a(x)}^{b(t)}v^{1-p^{\prime}}(y)dy\right)^{\frac{s-1}{p}}<\infty
		\end{equation}
holds with $s\in(1,p)$,
		where supremum is taken over all $x$ and $t$ such that $\eqref{eq:13}$ holds. Moreover, the best constant $C$ in $\eqref{eq:3}$ satisfies
		\begin{equation}\label{eq:5}
		\sup_{s\in (1,p)}B(s)\left(\frac{\left(\frac{p}{p-s}\right)^{p}}{\left(\frac{p}{p-s}\right)^{p}+\frac{1}{s-1}}\right)^{\frac{1}{p}}\leq C\leq2\inf_{s\in (1,p)}B(s)\left(\frac{p-1}{p-s}\right)^{\frac{1}{p^{\prime}}}.
		\end{equation}
		Specially, when $s=p$, Heinig and Sinnamon \cite{HS}  characterized the following sufficient and necessary condition:
		\begin{equation}\label{eq:6}
		A=\lim\limits_{s\rightarrow p}B(s)<\infty,
		\end{equation}
 and the best constant $C$ in $\eqref{eq:3}$ satisfies
		\begin{equation*}
		A\leq C\leq 2\left(1+\frac{q}{p^{\prime}}\right)^{\frac{1}{q}}\left(1+\frac{p^{\prime}}{q}\right)^{\frac{1}{p^{\prime}}}A,
		\end{equation*}
where $B(s)$ is as in (\ref{eq:4}).

		In this paper  we wish to characterize the sufficient conditions  ensure the boundedness of  the  operator $H_2$ in (\ref{eq:1}).

 Kufner and Persson \cite{KP} studied Pl\'{o}ya-Knnop type operator with variable limits  defined by
\begin{equation*}
(Gf)(x)=\exp\left(\frac{1}{b(x)-a(x)}\displaystyle\int_{a(x)}^{b(x)}\ln f(t)dt\right).
\end{equation*}
They obtained following inequality
\begin{equation}\label{eq:7}
\|(Gf)\|_{L^{q}(\mathbb{R}_{+},u)}\leq C\|f\|_{L^{p}(\mathbb{R}_{+},v)},
\end{equation}
 for  $0<p\leq q<\infty$ with the following weight condition
\begin{equation}\label{eq:8}
A_{PP}=\sup_{t>0}\left(\int_{\sigma^{-1}(a(x))}^{\sigma^{-1}(b(x))}w(x)dx\right)^{\frac{1}{q}}(b(t)-a(t))^{-\frac{1}{p}}<\infty
\end{equation}
holds,
where
$$\sigma(t)=\frac{a(t)+b(t)}{2},$$
and
\begin{equation}\label{eq:9}
w(x)=\exp\left(\frac{1}{b(x)-a(x)}\int_{a(x)}^{b(x)}\ln\frac{1}{v(t)}dt\right)^{\frac{q}{p}}u(x).
\end{equation}
Ushakova \cite{U} characterized a completely different condition rather than $\eqref{eq:8}$. That is, for any $s>1$,
$$D(s)=\sup\left(\int_{t}^{x}w(y)(b(y)-a(y))^{\frac{-qs}{p}}dy\right)^{\frac{1}{q}}(b(t)-a(x))^{\frac{s-1}{p}}<\infty,$$
where the supremum is taken over all $x$ and $t$ such that $\eqref{eq:13}$ holds and $w(x)$ is defined by $\eqref{eq:9}$.

	The corresponding two dimensional Pl\'{o}ya-Knnop type operator with variable limits is defined as
		\begin{equation}\label{eq:10}
		(G_{2}f)(x_{1},x_{2})=\exp\left(\frac{1}{\prod_{i=1}^2(b_{i}(x_{i})-a_{i}(x_{i}))}\int_{a_{1}(x_{1})}^{b_{1}(x_{1})}\int_{a_{2}(x_{2})}^{b_{2}(x_{2})}\ln f(t_{1},t_{2})dt_{1}dt_{2}\right).
		\end{equation}
	When $a_{1}(x_{1})=a_{2}(x_{2})=0$, $b_{1}(x_{1})=x_{1}$, $b_{2}(x_{2})=x_{2}$, Wedestig \cite{W} gave the characterization of sufficient and necessity conditions for the boundedness of the Pl\'{o}ya-Knnop operator $G_{2}$. That is, the following weight condition holds:
	$$D_{w}(s_{1},s_{2})=\sup_{{y_{1},y_{2}\in(0,\infty)}}y_{1}^{\frac{s_{1}-1}{p}}y_{2}^{\frac{s_{2}-1}{p}}\left(\int_{y_{1}}^{\infty}
\int_{y_{2}}^{\infty}x_{1}^{-\frac{s_{1}q}{p}}x_{2}^{-\frac{s_{2}q}{p}}w(x_{1},x_{2})dx_{1}dx_{2}\right)^{\frac{1}{q}}<\infty,$$
	where $s_{1},s_{2}>1$ and $w(x_{1},x_{2})$ is defined as
\begin{equation}\label{eq_w} w(x_{1},x_{2})=\exp\left(\frac{1}{\prod_{i=1}^{2}(b_{i}(x_{i})-a_{i}(x_{i}))}\int_{a_{1}(x_{1})}^{b_{1}(x_{1})}\int_{a_{2}(x_{2})}^{b_{2}(x_{2})}
\ln\frac{1}{v(t_{1},t_{2})}dt_{1}dt_{2}\right)^{\frac{q}{p}}u(x_{1},x_{2}).
	\end{equation}
	
	In this paper we wish to obtain sufficient and necessity conditions for the boundedness of the Pl\'{o}ya-Knnop operator $G_{2}$ with general $a_{i}(x_{i})$ and $b_{i}(x_{i})$ satisfy condition $\eqref{eq:2}$, $i=1,2$.
	
		{\centering \section{A scale of weighted characterization for  Hardy operator with viriable limits}}
		
		Our main result in this section reads:
		\begin{thm}\label{main:1}
			Let $1<p\leq q<\infty$, $1<s_{i}<p$ and let $u$ be a weight function on $\mathbb{R}_{+}^{2}:=\mathbb{R}_{+}\times \mathbb{R}_{+}$ and $v_{i}$ be weight function functions on $\mathbb{R}_{+}$ with $i=1,2$. Then the inequality
			\begin{equation}\label{eq:11}
			\|(H_{2}f)\|_{L^{q}(\mathbb{R}_{+}^{2},u)}\leq C_{p,q}\|f\|_{L^{p}(\mathbb{R}_{+}^{2},v_{1}v_{2})}
			\end{equation}
			holds for some finite constant $C_{p,q}$ and for all measurable functions $f\geq 0$ if and only if
			\begin{equation}\label{eq:12}
			{
				\begin{split}
				B(s_{1},s_{2})&:=\sup\left(\int_{t_{1}}^{x_{1}}\int_{t_{2}}^{x_{2}}u(y_{1},y_{2})\prod_{i=1}^{2}\left(V_{i}(b_{i}(y_{i}))-V_{i}(a_{i}(y_{i}))\right)^{\frac{q(p-s_{i})}{p}}dy_{1}dy_{2}\right)^{\frac{1}{q}}\\
				&\qquad\times\prod_{i=1}^{2}\left(V_{i}(b_{i}(t_{i}))-V_{i}(a_{i}(x_{i}))\right)^{\frac{(s_{i}-1)}{p}}<\infty,
				\end{split}
			}
			\end{equation}
			where
\begin{equation}\label{eq:V-F}
V_{i}(t_{i})=\int_{0}^{t_{i}}v_{i}^{1-p^{\prime}}(\eta)d\eta
\end{equation}
 and supremum is taken over all $x_{i}$ and $t_{i}$ such that
			\begin{equation}\label{eq:13}
			0<t_{i}<x_{i}<\infty\qquad   \text{and} \qquad a_{i}(x_{i})<b_{i}(t_{i}), \quad i=1,2.
			\end{equation}
			Furthermore, if $C_{p,q}$ is the best possible constant  in $\eqref{eq:11}$, then
			\begin{equation}\label{eq:14}
			\sup_{s_{1},s_{2}\in (1,p)}B(s_{1},s_{2})\prod_{i=1}^{2}\left(\frac{\left(\frac{p}{p-s_{i}}\right)^{p}}{\left(\frac{p}{p-s_{i}}\right)^{p}+\frac{1}{s_{i}-1}}\right)^{\frac{1}{p}}\leq C_{p,q}\leq4\inf_{s_{1},s_{2}\in (1,p)}B(s_{1},s_{2})\prod_{i=1}^{2}\left(\frac{p-1}{p-s_{i}}\right)^{\frac{1}{p^{\prime}}}.
			\end{equation}
		\end{thm}
		
To prove Theorem $\ref{main:1}$, we need the following result of Wdedestig \cite{W}:

\noindent {\bf Lemma A.}
			Let $1<p\leq q<\infty$,  $1<s_{i}<p$ and let $u$ be a weight function on $\mathbb{R}_{+}^{2}$ and $v_{i}$ be weight function functions on $\mathbb{R}_{+}$ with $i=1,2$. Then for $0\leq c_{i}<d_{i}\leq\infty$ the inequality
			\begin{equation}\label{eq:15}
			{
				\begin{split}
				&\left(\int_{c_{1}}^{d_{1}}\int_{c_{2}}^{d_{2}}\left(\int_{c_{1}}^{x_{1}}\int_{c_{2}}^{x_{2}}f(t_{1},t_{2})dt_{1}dt_{2}\right)^{q}u(x_{1},x_{2})dx_{1}dx_{2}\right)^{\frac{1}{q}}\\
				&\qquad\leq C	\left(\int_{c_{1}}^{d_{1}}\int_{c_{2}}^{d_{2}}f^{p}(x_{1},x_{2})v_{1}(x_{1})v_{2}(x_{2})dx_{1}dx_{2}\right)^{\frac{1}{p}}
				\end{split}
				}
			\end{equation}
		holds for some finite constant $C$ and for all measurable functions $f\geq 0$ if and only if
		
		\begin{equation}\label{eq:16}
			{
				\begin{split}
				A_{W}(s_{1},s_{2})&=\sup_{\substack{c_{1}\leq t_{1}\leq d_{1}\\c_{2}\leq t_{2}\leq d_{2}}}\left(\int_{t_{1}}^{d_{1}}\int_{t_{2}}^{d_{2}}u(x_{1},x_{2})\prod_{i=1}^{2}\left(V_{i}(x_{i})-V_{i}(c_{i})\right)^{\frac{q(p-s_{i})}{p}}dx_{1}dx_{2}\right)^{\frac{1}{q}}\\
				&\qquad\times\prod_{i=1}^{2}\left(V_{i}(t_{i})-V_{i}(c_{i})\right)^{\frac{(s_{i}-1)}{p}}<\infty,
				\end{split}
			}
			\end{equation}
		where the function $V_{i}(x_{i})$ is defined by $\eqref{eq:V-F}$ with $i=1,2$.		
		Furthermore, if $C$ is the best possible constant  in $\eqref{eq:15}$, then
		\begin{equation}\label{eq:17}
		\sup_{s_{1},s_{2}\in (1,p)}A_{W}(s_{1},s_{2})\prod_{i=1}^{2}\left(\frac{\left(\frac{p}{p-s_{i}}\right)^{p}}{\left(\frac{p}{p-s_{i}}\right)^{p}+\frac{1}{s_{i}-1}}\right)^{\frac{1}{p}}\leq C\leq\inf_{s_{1},s_{2}\in (1,p)}A_{W}(s_{1},s_{2})\prod_{i=1}^{2}\left(\frac{p-1}{p-s_{i}}\right)^{\frac{1}{p^{\prime}}}.
		\end{equation}

			\begin{lem}\label{lem:2}
				Let $1<p\leq q<\infty$,  $1<s_{i}<p$ and let $u$ be a weight function on $\mathbb{R}_{+}^{2}$ and $v_{i}$ be weight function functions on $\mathbb{R}_{+}$ with $i=1,2$. Then for $0\leq c_{i}<d_{i}\leq\infty$ the inequality
				\begin{equation}\label{eq:18}
				{
					\begin{split}
					&\left(\int_{c_{1}}^{d_{1}}\int_{c_{2}}^{d_{2}}\left(\int_{c_{1}}^{x_{1}}\int_{x_{2}}^{d_{2}}f(t_{1},t_{2})dt_{1}dt_{2}\right)^{q}u(x_{1},x_{2})dx_{1}dx_{2}\right)^{\frac{1}{q}}\\
					&\qquad\leq C	\left(\int_{c_{1}}^{d_{1}}\int_{c_{2}}^{d_{2}}f^{p}(x_{1},x_{2})v_{1}(x_{1})v_{2}(x_{2})dx_{1}dx_{2}\right)^{\frac{1}{p}}
					\end{split}
				}
				\end{equation}
				holds for some finite constant $C$ and for all measurable functions $f\geq 0$ if and only if
				
				\begin{equation}\label{eq:19}
				{
					\begin{split}
					&	A_{W}^{*}(s_{1},s_{2})=\sup_{\substack{c_{1}\leq t_{1}\leq d_{1}\\c_{2}\leq t_{2}\leq d_{2}}}(V_{1}(t_{1})-V_{1}(c_{1}))^{\frac{(s_{1}-1)}{p}}(V_{2}(d_{2})-V_{2}(t_{2}))^{\frac{(s_{2}-1)}{p}}\\
					&\times
					\left(\int_{t_{1}}^{d_{1}}\int_{c_{2}}^{t_{2}}u(x_{1},x_{2})(V_{1}(x_{1})-V_{1}(c_{1}))^{\frac{q(p-s_{1})}{p}}(V_{2}(d_{2})-V_{2}(x_{2}))^{\frac{q(p-s_{2})}{p}}dx_{1}dx_{2}\right)^{\frac{1}{q}}<\infty,
					\end{split}
				}
				\end{equation}
				where the function $V_{i}(x_{i})$ is defined by $\eqref{eq:V-F}$ with $i=1,2$.				
				Furthermore, if $C$ is the best possible constant  in $\eqref{eq:18}$, then
				\begin{equation}\label{eq:20}
				\sup_{s_{1},s_{2}\in (1,p)}A_{W}^{*}(s_{1},s_{2})\prod_{i=1}^{2}\left(\frac{\left(\frac{p}{p-s_{i}}\right)^{p}}{\left(\frac{p}{p-s_{i}}\right)^{p}+\frac{1}{s_{i}-1}}\right)^{\frac{1}{p}}\leq C\leq\inf_{s_{1},s_{2}\in (1,p)}A_{W}^{*}(s_{1},s_{2})\prod_{i=1}^{2}\left(\frac{p-1}{p-s_{i}}\right)^{\frac{1}{p^{\prime}}}.
				\end{equation}
			\end{lem}
		
Proof. Let
$$g(x_{1},x_{2})=f^{p}(x_{1},x_{2})v_{1}(x_{1})v_{2}(x_{2})$$
 in $\eqref{eq:18}$. The inequality $\eqref{eq:18}$ is equivalent to the following inequality
			\begin{equation}\label{eq:21}
			{
				\begin{split}
				&\left(\int_{c_{1}}^{d_{1}}\int_{c_{2}}^{d_{2}}\left(\int_{c_{1}}^{x_{1}}\int_{x_{2}}^{d_{2}}g^{\frac{1}{p}}(t_{1},t_{2})v_{1}^{-\frac{1}{p}}(t_{1})
v_{2}^{-\frac{1}{p}}(t_{2})dt_{1}dt_{2}\right)^{q}u(x_{1},x_{2})dx_{1}dx_{2}\right)^{\frac{1}{q}}\\
				&\qquad\leq C	\left(\int_{c_{1}}^{d_{1}}\int_{c_{2}}^{d_{2}}g(x_{1},x_{2})dx_{1}dx_{2}\right)^{\frac{1}{p}}.
				\end{split}
			}
			\end{equation}
	Assume  $\eqref{eq:19}$	holds.
	$$\frac{d}{dt_{1}}V_{1}(t_{1})=v_{1}^{1-p^{\prime}}=v_{1}^{-\frac{p^{\prime}}{p}}(t_{1}),\quad \frac{d}{dt_{2}}V_{2}(t_{2})=v_{2}^{1-p^{\prime}}=v_{2}^{-\frac{p^{\prime}}{p}}(t_{2}).$$
It follows from  H\"{o}lder's inequality and Minkowski's inequality  that	
$$
{
	\begin{split}
	&\left(\int_{c_{1}}^{d_{1}}\int_{c_{2}}^{d_{2}}\left(\int_{c_{1}}^{x_{1}}\int_{x_{2}}^{d_{2}}g^{\frac{1}{p}}(t_{1},t_{2})v_{1}^{-\frac{1}{p}}(t_{1})v_{2}^{
	-\frac{1}{p}}(t_{2})dt_{1}dt_{2}\right)^{q}u(x_{1},x_{2})dx_{1}dx_{2}\right)^{\frac{1}{q}}\\
	&\quad=\left(\int_{c_{1}}^{d_{1}}\int_{c_{2}}^{d_{2}}\left(\int_{c_{1}}^{x_{1}}\int_{x_{2}}^{d_{2}}g^{\frac{1}{p}}(t_{1},t_{2})\left((V_{1}(t_{1})-V_{1}(c_{1})\right)^{\frac{s_{1}-1}{p}}
\left(V_{2}(d_{2})-V_{2}(t_{2})\right)^{\frac{s_{2}-1}{p}}\right.\right.\\
	&\qquad\times\left.\left.\left((V_{1}(t_{1})-V_{1}(c_{1})\right)^{-\frac{s_{1}-1}{p}}v_{1}^{-\frac{1}{p}}(t_{1})\left(V_{2}(d_{2})-V_{2}(t_{2})\right)^{\frac{s_{2}-1}{p}}v_{2}^{
		-\frac{1}{p}}(t_{2})dt_{1}dt_{2}\right)^{q}u(x_{1},x_{2})dx_{1}dx_{2}\right)^{\frac{1}{q}}\\
	&\quad=\left(\int_{c_{1}}^{d_{1}}\int_{c_{2}}^{d_{2}}\left(\int_{c_{1}}^{x_{1}}\int_{x_{2}}^{d_{2}}g(t_{1},t_{2})V_{1,c_{1}}^{s_{1}-1}(t_{1})V_{2,d_{2}}^{s_{2}-1}(t_{2})dt_{1}dt_{2}\right)^{\frac{q}{p}}\right.\\
		\end{split}
	}
	$$	
	$$
	{
		\begin{split}
	&\qquad\times\left.\left(\int_{c_{1}}^{x_{1}}V_{1,c_{1}}^{-\frac{(s_{1}-1)p^{\prime}}{p}}(t_{1})v_{1}^{-\frac{p^{\prime}}{p}}(t_{1})dt_{1}\right)^{\frac{q}{p^{\prime}}}
\left(\int_{x_{2}}^{d_{2}}V_{2,d_{2}}^{-\frac{(s_{2}-1)p^{\prime}}{p}}(t_{2})v_{2}^{-\frac{p^{\prime}}{p}}(t_{2})dt_{2}\right)^{\frac{q}{p^{\prime}}}u(x_{1},x_{2})dx_{1}dx_{2}\right)^{\frac{1}{q}}\\
	&\quad=\left(\frac{p-1}{p-s_{1}}\right)^{\frac{1}{p^{\prime}}}\left(\frac{p-1}{p-s_{2}}\right)^{\frac{1}{p^{\prime}}}\left(\int_{c_{1}}^{d_{1}}\int_{c_{2}}^{d_{2}}
\left(\int_{c_{1}}^{x_{1}}\int_{x_{2}}^{d_{2}}g(t_{1},t_{2})V_{1,c_{1}}^{s_{1}-1}(t_{1})V_{2,d_{2}}^{s_{2}-1}(t_{2})dt_{1}dt_{2}\right)^{\frac{q}{p}}\right.\\
	&\qquad\times\left.V_{1,c_{1}}^{\frac{q(p-s_{1})}{p}}(x_{1})V_{2,d_{2}}^{\frac{q(p-s_{2})}{p}}(x_{2})u(x_{1},x_{2})dx_{1}dx_{2}\right)^{\frac{1}{q}}\\
	&\quad=\left(\frac{p-1}{p-s_{1}}\right)^{\frac{1}{p^{\prime}}}\left(\frac{p-1}{p-s_{2}}\right)^{\frac{1}{p^{\prime}}}
\left(\int_{c_{1}}^{d_{1}}\int_{c_{2}}^{d_{2}}g(t_{1},t_{2})V_{1,c_{1}}^{s_{1}-1}(t_{1})V_{2,d_{2}}^{s_{2}-1}(t_{2})dt_{1}dt_{2}\right.\\
	&\qquad\times\left.\left(\int_{x_{1}}^{d_{1}}\int_{c_{2}}^{x_{2}}V_{1,c_{1}}^{\frac{q(p-s_{1})}{p}}(x_{1})V_{2,d_{2}}^{\frac{q(p-s_{2})}{p}}(x_{2})u(x_{1},x_{2})dx_{1}dx_{2}\right)^{\frac{p}{q}}dt_{1}dt_{2}\right)^{\frac{1}{p}}\\
	&\quad\leq\left(\frac{p-1}{p-s_{1}}\right)^{\frac{1}{p^{\prime}}}\left(\frac{p-1}{p-s_{2}}\right)^{\frac{1}{p^{\prime}}}A_{W}^{*}(s_{1},s_{2})
\left(\int_{c_{1}}^{d_{1}}\int_{c_{2}}^{d_{2}}g(t_{1},t_{2})dt_{1}dt_{2}\right)^{\frac{1}{p}},
	\end{split}
	}
$$
where $V_{1,c_{1}}(t_{1})=v_{1}(t_{1})-V_{1}(c_{1})$ and $V_{2,d_{2}}(t_{2})=V_{2}(d_{2})-V_{2}(t_{2})$.

 Hence, $\eqref{eq:21}$ and $\eqref{eq:18}$ holds with a constant satisfying the right hand side inequality in $\eqref{eq:20}$.

 Now we assume that $\eqref{eq:18}$ and, thus, $\eqref{eq:21}$ holds. Consider the test function
$$
{
	\begin{split}
		&g(x_{1},x_{2})=\left(\frac{p}{p-s_{1}}\right)^{p}\left(\frac{p}{p-s_{2}}\right)^{p}\\
		&\times\left(V_{1}(t_{1})-V_{1}(c_{1})\right)^{-s_{1}}v_{1}^{1-p^{\prime}}(x_{1})\left(V_{2}(d_{2})-V_{2}(t_{2})\right)^{-s_{2}}v_{2}^{1-p^{\prime}}(x_{2})\chi_{(c_{1},t_{1})}(x_{1})\chi_{(t_{2},d_{2})}(x_{2})\\
		&+\left(\frac{p}{p-s_{1}}\right)^{p}\left(V_{1}(t_{1})-V_{1}(c_{1})\right)^{-s_{1}}v_{1}^{1-p^{\prime}}(x_{1})
\left(V_{2}(d_{2})-V_{2}(x_{2})\right)^{-s_{2}}v_{2}^{1-p^{\prime}}(x_{2})\chi_{(c_{1},t_{1})}(x_{1})\chi_{(c_{2},t_{2})}(x_{2})\\
		&+\left(\frac{p}{p-s_{2}}\right)^{p}\left(V_{1}(x_{1})-V_{1}(c_{1})\right)^{-s_{1}}v_{1}^{1-p^{\prime}}(x_{1})
\left(V_{2}(d_{2})-V_{2}(t_{2})\right)^{-s_{2}}v_{2}^{1-p^{\prime}}(x_{2})\chi_{(t_{1},d_{1})}(x_{1})\chi_{(t_{2},d_{2})}(x_{2})\\
		&+\left(V_{1}(x_{1})-V_{1}(c_{1})\right)^{-s_{1}}v_{1}^{1-p^{\prime}}(x_{1})\left(V_{2}(d_{2})-V_{2}(x_{2})\right)^{-s_{2}}v_{2}^{1-p^{\prime}}(x_{2})\chi_{(t_{1},d_{1})}(x_{1})\chi_{(c_{2},t_{2})}(x_{2}),
	\end{split}
	}
$$
where $t_{1}$, $t_{2}$ are fix numbers and $c_{i}\leq t_{i}\leq d_{i}$ with $i=1,2$. Then we show that the integral on right hand side of $\eqref{eq:21}$
$$
{
	\begin{split}
	&\qquad\qquad\qquad\qquad\qquad\qquad\left(\int_{c_{1}}^{d_{1}}\int_{c_{2}}^{d_{2}}g(x_{1},x_{2})dx_{1}dx_{2}\right)^{\frac{1}{p}}=\\
	&\quad\left(\prod_{i=1}^{2}\left(\frac{p}{p-s_{i}}\right)^{p}\int_{c_{1}}^{t_{1}}\left(V_{1}(t_{1})-V_{1}(c_{1})\right)^{-s_{1}}v_{1}^{1-p^{\prime}}(x_{1})dx_{1}\int_{t_{2}}^{d_{2}}
\left(V_{2}(d_{1})-V_{2}(t_{2})\right)^{-s_{2}}v_{2}^{1-p^{\prime}}(x_{2})dx_{2}\right.\\
	&\quad+\int_{c_{1}}^{t_{1}}\left(\frac{p}{p-s_{1}}\right)^{p}\left(V_{1}(t_{1})-V_{1}(c_{1})\right)^{-s_{1}}v_{1}^{1-p^{\prime}}(x_{1})dx_{1}
\int_{c_{2}}^{t_{2}}\left(V_{2}(d_{2})-V_{2}(x_{2})\right)^{-s_{2}}v_{2}^{1-p^{\prime}}(x_{2})dx_{2}\\
		\end{split}
	}
	$$	
	$$
	{
		\begin{split}
	&\quad+\int_{t_{1}}^{d_{1}}\left(V_{1}(d_{1})-V_{1}(x_{1})\right)^{-s_{1}}v_{1}^{1-p^{\prime}}(x_{1})dx_{1}\int_{t_{2}}^{d_{2}}\left(\frac{p}{p-s_{2}}\right)^{p}
\left(V_{2}(t_{2})-V_{2}(c_{2})\right)^{-s_{2}}v_{2}^{1-p^{\prime}}(x_{2})dx_{2}\\
	&\quad+\left.\int_{t_{1}}^{d_{1}}\left(V_{1}(d_{1})-V_{1}(x_{1})\right)^{-s_{1}}v_{1}^{1-p^{\prime}}(x_{1})dx_{1}
\int_{c_{2}}^{t_{2}}\left(V_{2}(d_{2})-V_{1}(x_{2})\right)^{-s_{2}}v_{2}^{1-p^{\prime}}(x_{2})dx_{2}\right)^{\frac{1}{p}}\\
	&\leq\left(\left(\frac{p}{p-s_{1}}\right)^{p}+\frac{1}{s_{1}-1}\right)^{\frac{1}{p}}
\left(\left(\frac{p}{p-s_{2}}\right)^{p}+\frac{1}{s_{2}-1}\right)^{\frac{1}{p}}\left(V_{1}(t_{1})-V_{1}(c_{1})\right)^{\frac{1-s_{1}}{p}}\left(V_{2}(d_{2})-V_{2}(t_{2})\right)^{\frac{1-s_{2}}{p}}.
	\end{split}
	}
$$
Moreover, the left hand side of $\eqref{eq:21}$ is greater than
$$
{
	\begin{split}
	&\left(\int_{t_{1}}^{d_{1}}\int_{c_{1}}^{t_{1}}\left[\prod_{i=1}^{2}\frac{p}{p-s_{i}}\int_{c_{1}}^{t_{1}}\int_{t_{2}}^{d_{2}}
\left(V_{1}(t_{1})-V_{1}(c_{1})\right)^{-\frac{s_{1}}{p}}\left(V_{2}(d_{2})-V_{2}(t_{2})\right)^{-\frac{s_{2}}{p}}v_{1}^{1-p^{\prime}}(y_{1})v_{2}^{1-p^{\prime}}(y_{2})dy_{1}dy_{2}\right.\right.\\
	&\quad+\left(\int_{c_{1}}^{t_{1}}\frac{p}{p-s_{1}}\left(V_{1}(t_{1})-V_{1}(c_{1})\right)^{-\frac{s_{1}}{p}}v_{1}^{1-p^{\prime}}(y_{1})dy_{1}\right)\left(\int_{x_{2}}^{t_{2}}
\left(V_{2}(y_{2})-V_{2}(t_{2})\right)^{-\frac{s_{2}}{p}}v_{2}^{1-p^{\prime}}(y_{2})dy_{2}\right)\\
	&\quad+\left(\int_{t_{1}}^{x_{1}}\left(V_{1}(y_{1})-V_{1}(c_{1})\right)^{-\frac{s_{1}}{p}}v_{1}^{1-p^{\prime}}(y_{1})dy_{1}\right)
\left(\int_{t_{2}}^{d_{2}}\frac{p}{p-s_{2}}\left(V_{2}(d_{2})-V_{2}(t_{2})\right)^{-\frac{s_{2}}{p}}v_{2}^{1-p^{\prime}}(y_{2})dy_{2}\right)\\
	&\quad+\left.\left(\int_{t_{1}}^{x_{1}}\left(V_{1}(y_{1})-V_{1}(c_{1})\right)^{-\frac{s_{1}}{p}}v_{1}^{1-p^{\prime}}(y_{1})dy_{1}\right)
\left(\int_{x_{2}}^{t_{2}}\left(V_{2}(d_{2})-V_{2}(t_{2})\right)^{-\frac{s_{2}}{p}}v_{2}^{1-p^{\prime}}(y_{2})dy_{2}\right)\right]^{q}\\
	&\quad\times\left. u(x_{1},x_{2})dx_{1}dx_{2}\right)^{\frac{1}{q}}\\
	&=\frac{p}{p-s_{1}}\frac{p}{p-s_{2}}\left(\int_{t_{1}}^{d_{1}}\int_{c_{2}}^{t_{2}}u(x_{1},x_{2})\left(V_{1}(x_{1})-V_{1}(c_{1})\right)^{\frac{q(p-s_{1})}{p}}
\left(V_{2}(d_{1})-V_{2}(x_{2})\right)^{\frac{q(p-s_{2})}{p}}dx_{1}dx_{2}\right)^{\frac{1}{q}}.
	\end{split}
	}
$$
Hence, $\eqref{eq:21}$ implies that
\begin{equation}\label{eq:as}
{
	\begin{split}
	&\frac{p}{p-s_{1}}\frac{p}{p-s_{2}}\left(\int_{t_{1}}^{d_{1}}\int_{c_{2}}^{t_{2}}u(x_{1},x_{2})\left(V_{1}(x_{1})-V_{1}(c_{1})\right)^{\frac{q(p-s_{1})}{p}}
\left(V_{2}(d_{1})-V_{2}(x_{2})\right)^{\frac{q(p-s_{2})}{p}}dx_{1}dx_{2}\right)^{\frac{1}{q}}\\
	&\leq C\prod_{i=1}^{2}\left(\left(\frac{p}{p-s_{i}}\right)^{p}+\frac{1}{s_{i}-1}\right)^{\frac{1}{p}}\left(V_{1}(t_{1})-V_{1}(c_{1})\right)^{\frac{1-s_{1}}{p}}\left(V_{2}(d_{2})-V_{2}(t_{2})\right)^{\frac{1-s_{2}}{p}}.
	\end{split}
	}
\end{equation}
Above $\eqref{eq:as}$ inequality is equivalent to
$$
{
	\begin{split}
&\left(\int_{t_{1}}^{d_{1}}\int_{c_{2}}^{t_{2}}u(x_{1},x_{2})\left(V_{1}(x_{1})-V_{1}(c_{1})\right)^{\frac{q(p-s_{1})}{p}}\left(V_{2}(d_{1})-V_{2}(x_{2})\right)^{\frac{q(p-s_{2})}{p}}dx_{1}dx_{2}\right)^{\frac{1}{q}}\\
&\times\left(V_{1}(t_{1})-V_{1}(c_{1})\right)^{\frac{1-s_{1}}{p}}\left(V_{2}(d_{2})-V_{2}(t_{2})\right)^{\frac{1-s_{2}}{p}}\prod_{i=1}^{2}
\left(\frac{\left(\frac{p}{p-s_{i}}\right)^{p}}{\left(\frac{p}{p-s_{i}}\right)^{p}+\frac{1}{s_{i}-1}}\right)^{\frac{1}{p}}\leq C.
	\end{split}
	}	
	$$
We obtain that $\eqref{eq:19}$ and the left hand side of the estimate of $\eqref{eq:20}$ hold. $\hfill \Box $

By applying a standard duality, argument Kufner and Persson \cite{KP} we have that
		\begin{lem}\label{lem:3}
			Let $1<p\leq q<\infty$,  $1<s_{i}<p$ and let $u$ be a weight function on $\mathbb{R}_{+}^{2}$ and $v_{i}$ be weight function functions on $\mathbb{R}_{+}$ with $i=1,2$. Then for $0\leq c_{i}<d_{i}\leq\infty$ the inequality
			\begin{equation}\label{eq:22}
			{
				\begin{split}
				&\left(\int_{c_{1}}^{d_{1}}\int_{c_{2}}^{d_{2}}\left(\int_{x_{1}}^{d_{1}}\int_{x_{2}}^{d_{2}}f(t_{1},t_{2})dt_{1}dt_{2}\right)^{q}u(x_{1},x_{2})dx_{1}dx_{2}\right)^{\frac{1}{q}}\\
				&\qquad\qquad\qquad\qquad\qquad\qquad\qquad\leq C	\left(\int_{c_{1}}^{d_{1}}\int_{c_{2}}^{d_{2}}f^{p}(x_{1},x_{2})v_{1}(x_{1})v_{2}(x_{2})dx_{1}dx_{2}\right)^{\frac{1}{p}}
				\end{split}
			}
			\end{equation}
			holds for some finite constant $C$ and for all measurable functions $f\geq 0$ if and only if
			
			\begin{equation}\label{eq:23}
			{
				\begin{split}
					\tilde{A}_{W}(s_{1},s_{2})&=\sup_{\substack{c_{1}\leq t_{1}\leq d_{1}\\c_{2}\leq t_{2}\leq d_{2}}}\left(\int_{c_{1}}^{t_{1}}\int_{c_{2}}^{t_{2}}u(x_{1},x_{2})\prod_{i=1}^{2}\left(V_{i}(d_{i})-V_{i}(d_{i})\right)^{\frac{q(p-s_{i})}{p}}dx_{1}dx_{2}\right)^{\frac{1}{q}}\\
				&\qquad\times\prod_{i=1}^{2}\left(V_{i}(d_{i})-V_{i}(t_{i})\right)^{\frac{(s_{i}-1)}{p}}<\infty,
				\end{split}
			}
			\end{equation}
			where the function $V_{i}(x_{i})$ is defined by in $\eqref{eq:V-F}$ with $i=1,2$.
			
			Furthermore, if $C$ is the best possible constant  in $\eqref{eq:22}$, then
			\begin{equation}\label{eq:24}
			\sup_{s_{1},s_{2}\in (1,p)}\tilde{A}_{W}(s_{1},s_{2})\prod_{i=1}^{2}\left(\frac{\left(\frac{p}{p-s_{i}}\right)^{p}}{\left(\frac{p}{p-s_{i}}\right)^{p}+\frac{1}{s_{i}-1}}\right)^{\frac{1}{p}}\leq C\leq\inf_{s_{1},s_{2}\in (1,p)}\tilde{A}_{W}(s_{1},s_{2})\prod_{i=1}^{2}\left(\frac{p-1}{p-s_{i}}\right)^{\frac{1}{p^{\prime}}}.
			\end{equation}	
		\end{lem}

		\begin{lem}\label{lem:4}
			Let $1<p\leq q<\infty$,  $1<s_{i}<p$ and let $u$ be a weight function on $\mathbb{R}_{+}^{2}$ and $v_{i}$ be weight function functions on $\mathbb{R}_{+}$ with $i=1,2$. Then for $0\leq c_{i}<d_{i}\leq\infty$ the inequality
			\begin{equation}\label{eq:25}
			{
				\begin{split}
				&\left(\int_{c_{1}}^{d_{1}}\int_{c_{2}}^{d_{2}}\left(\int_{x_{1}}^{d_{1}}\int_{c_{2}}^{x_{2}}f(t_{1},t_{2})dt_{1}dt_{2}\right)^{q}u(x_{1},x_{2})dx_{1}dx_{2}\right)^{\frac{1}{q}}\\
				&\qquad\qquad\qquad\qquad\qquad\qquad\qquad\leq C	\left(\int_{c_{1}}^{d_{1}}\int_{c_{2}}^{d_{2}}f^{p}(x_{1},x_{2})v_{1}(x_{1})v_{2}(x_{2})dx_{1}dx_{2}\right)^{\frac{1}{p}}
				\end{split}
			}
			\end{equation}
			holds for some finite constant $C$ and for all measurable functions $f\geq 0$ if and only if
			
			\begin{equation}\label{eq:26}
			{
				\begin{split}
				&	\tilde{A}_{W}^{*}(s_{1},s_{2})=\sup_{\substack{c_{1}\leq t_{1}\leq d_{1}\\c_{2}\leq t_{2}\leq d_{2}}}(V_{1}(d_{1})-V_{1}(t_{1}))^{\frac{(s_{1}-1)}{p}}(V_{2}(t_{2})-V_{2}(c_{2}))^{\frac{(s_{2}-1)}{p}}\\
				&\times
				\left(\int_{c_{1}}^{t_{1}}\int_{t_{2}}^{d_{2}}u(x_{1},x_{2})(V_{1}(d_{1})-V_{1}(x_{1}))^{\frac{q(p-s_{1})}{p}}(V_{2}(x_{2})-V_{2}(c_{2}))^{\frac{q(p-s_{2})}{p}}dx_{1}dx_{2}\right)^{\frac{1}{q}}<\infty,
				\end{split}
			}
			\end{equation}
			where the function $V_{i}(x_{i})$ is defined by in $\eqref{eq:V-F}$  with $i=1,2$.			
			Furthermore, if $C$ is the best possible constant  in $\eqref{eq:25}$, then we have
			\begin{equation}\label{eq:27}
			\sup_{s_{1},s_{2}\in (1,p)}\tilde{A}_{W}^{*}(s_{1},s_{2})\prod_{i=1}^{2}\left(\frac{\left(\frac{p}{p-s_{i}}\right)^{p}}{\left(\frac{p}{p-s_{i}}\right)^{p}+\frac{1}{s_{i}-1}}\right)^{\frac{1}{p}}\leq C\leq\inf_{s_{1},s_{2}\in (1,p)}\tilde{A}_{W}^{*}(s_{1},s_{2})\prod_{i=1}^{2}\left(\frac{p-1}{p-s_{i}}\right)^{\frac{1}{p^{\prime}}}.
			\end{equation}
		\end{lem}
	
	Next we begin to prove Theorem\ref{main:1}. We first prove the necessity. Assume that $\eqref{eq:11}$ holds.  We  choose the test function
	$$
	{
		\begin{split}
		f(y_{1},y_{2})&=\prod_{i=1}^{2}\frac{p}{p-s_{i}}\big(V_{i}(b_{i}(t_{i})-V_{i}(a_{i}(z_{i}))\big)^{-\frac{s_{i}}{p}}v_{i}^{1-p^{\prime}}(y_{i})\chi_{(a_{i}(z_{i}),b_{i}(t_{i}))}(y_{i})\\
		&+\frac{p}{p-s_{1}}\left(V_{1}(b_{1}(t_{1})-V_{1}(a_{1}(z_{1}))\right)^{-\frac{s_{1}}{p}}v_{1}^{1-p^{\prime}}(y_{1})\left(V_{2}(b_{2}(y_{2})-V_{2}(a_{2}(z_{2}))\right)^{-\frac{s_{2}}{p}}v_{2}^{1-p^{\prime}}(y_{2})\\
		&\times\chi_{(a_{1}(z_{1}),b_{1}(t_{1}))}(y_{1})\chi_{(b_{2}(t_{2}),b_{2}(z_{2}))}(y_{2})+\frac{p}{p-s_{2}}\left(V_{1}(b_{1}(y_{1})-V_{1}(a_{1}(z_{1}))\right)^{-\frac{s_{1}}{p}}v_{1}^{1-p^{\prime}}(y_{1})\\
		&\times\left(V_{2}(b_{2}(t_{2})-V_{2}(a_{2}(z_{2}))\right)^{-\frac{s_{2}}{p}}v_{2}^{1-p^{\prime}}(y_{2})\chi_{(b_{1}(t_{1}),b_{1}(z_{1}))}(y_{1})\chi_{(a_{2}(z_{2}),b_{2}(t_{2}))}(y_{2})\\
		&+\prod_{i=1}^{2}\left(V_{i}(b_{i}(y_{i})-V_{i}(a_{i}(z_{i}))\right)^{-\frac{s_{i}}{p}}v_{i}^{1-p^{\prime}}(y_{i})\chi_{(b_{i}(t_{i}),b_{i}(z_{i}))}(y_{i}),
		\end{split}
		}
	$$
	where $t_{i}$ and $x_{i}$ are fixed  numbers and  $0<t_{i}<x_{i}<\infty$ such that $a_{i}(x_{i})<b_{i}(t_{i})$ with $i=1,2$.

 It follows that  $$a_{i}(z_{i})\leq a_{i}(x_{i})<b_{i}(t_{i})<b_{i}(z_{i})\leq b_{i}(x_{i}),$$
 when   $t_{i}<z_{i}<x_{i}$ for $i=1,2$. Then  the right-hand side of $\eqref{eq:11}$ is equal to
	\begin{equation*}
	{
		\begin{split} &\prod_{i=1}^{2}\int_{a_{i}(z_{i})}^{b_{i}(t_{i})}\left(\frac{p}{p-s_{i}}\right)^{p}\left(V_{i}(b_{i}(t_{i}))-V_{i}(a_{i}(z_{i}))\right)^{-s_{i}}v_{i}^{1-p^{\prime}}(y_{i})dy_{i}+\\
		&\int_{a_{1}(z_{1})}^{b_{1}(t_{1})}\left(\frac{p}{p-s_{1}}\right)^{p}\left(V_{1}(b_{1}(t_{1}))-V_{1}(a_{1}(z_{1}))\right)^{-s_{1}}v_{1}^{1-p^{\prime}}(y_{1})dy_{1}
		\int_{b_{2}(t_{2})}^{b_{2}(z_{2})}\left(V_{2}(b_{2}(y_{2})-V_{2}(a_{2}(z_{2}))\right)^{-s_{2}}\\
		&\quad\times v_{2}^{1-p^{\prime}}(y_{2})dy_{2}+\int_{b_{1}(t_{1})}^{b_{1}(z_{1})}\left(V_{1}(b_{1}(y_{1}))-V_{1}(a_{1}(z_{1}))\right)^{-s_{1}}v_{1}^{1-p^{\prime}}(y_{1})dy_{1}
		\int_{a_{2}(z_{2})}^{b_{2}(t_{2})}\left(\frac{p}{p-s_{2}}\right)^{p}\\
		&\quad\times\left(V_{2}(b_{2}(t_{2})-V_{2}(a_{2}(z_{2}))\right)^{-\frac{s_{2}}{p}}v_{2}^{1-p^{\prime}}(y_{2})dy_{2}+\prod_{i=1}^{2}\int_{b_{i}(t_{i})}^{b_{i}(z_{i})}
\left(V_{i}(b_{i}(y_{i}))-V_{i}(a_{i}(z_{i}))\right)^{-s_{i}}v_{i}^{1-p^{\prime}}(y_{i})dy_{i}\\
		&\leq\prod_{i=1}^{2}\left(\left(\frac{p}{p-s_{i}}\right)^{p}+\frac{1}{s_{i}-1}\right)\left(V_{i}(b_{i}(t_{i}))-V_{i}(a_{i}(z_{i}))\right)^{1-s_{i}}.\\
\end{split}
	}
	\end{equation*}	
Notice that $a_{i}(z_{i})\leq a_{i}(x_{i}), i=1,2$.
We conclude that the integral on the right-hand side of $\eqref{eq:11}$ is less than or equal to
\begin{equation}\label{eq:2.18}
	{
		\begin{split}
\prod_{i=1}^{2}\left(\left(\frac{p}{p-s_{i}}\right)^{p}+\frac{1}{s_{i}-1}\right)\left(V_{i}(b_{i}(t_{i}))-V_{i}(a_{i}(x_{i}))\right)^{1-s_{i}}.
		\end{split}
	}
	\end{equation}

 On the other hand, for the left-hand side of $\eqref{eq:11}$, we have that
		\begin{equation}\label{eq:2.19}
		{
			\begin{split}
			&\left(\int_{0}^{\infty}\int_{0}^{\infty}\left(\int_{a_{1}(z_{1})}^{b_{1}(z_{1})}\int_{a_{2}(z_{2})}^{b_{2}(z_{2})}f(y_{1},y_{2})dy_{1}dy_{2}\right)^{q}u(z_{1},z_{2})dz_{1}dz_{2}\right)^{\frac{1}{q}}\\
			&\quad\geq\left(\int_{t_{1}}^{x_{1}}\int_{t_{2}}^{x_{2}}\left(\int_{a_{1}(z_{1})}^{b_{1}(z_{1})}\int_{a_{2}(z_{2})}^{b_{2}(z_{2})}f(y_{1},y_{2})dy_{1}dy_{2}\right)^{q}u(z_{1},z_{2})dz_{1}dz_{2}\right)^{\frac{1}{q}}\\
			&\quad=\left(\int_{t_{1}}^{x_{1}}\int_{t_{2}}^{x_{2}}\left(\prod_{i=1}^{2}\int_{a_{i}(z_{i})}^{b_{i}(t_{i})}\frac{p}{p-s_{i}}
\left(V_{i}(b_{i}(t_{i}))-V_{i}(a_{i}(z_{i}))\right)^{-\frac{s_{i}}{p}}v_{i}^{1-p^{\prime}}(y_{i})dy_{i}\right.\right.\\
			&\quad+\int_{a_{1}(z_{1})}^{b_{1}(t_{1})}\int_{b_{2}(t_{2})}^{b_{2}(z_{2})}\frac{p}{p-s_{1}}\left(V_{1}(b_{1}(t_{1}))-V_{1}(a_{1}(z_{1}))\right)^{-\frac{s_{1}}{p}}
\left(V_{2}(y_{2})-V_{1}(a_{2}(z_{2}))\right)^{-\frac{s_{2}}{p}}v_{1}^{1-p^{\prime}}(y_{1})\\
			&\qquad\times v_{2}^{1-p^{\prime}}(y_{2})dy_{1}dy_{2}\\
			&\quad+\int_{b_{1}(t_{1})}^{b_{1}(z_{1})}\int_{a_{2}(z_{2})}^{b_{2}(t_{2})}\frac{p}{p-s_{1}}\left(V_{1}(y_{1}))-V_{1}(a_{1}(z_{1}))\right)^{-\frac{s_{1}}{p}}
\left(V_{2}(b_{2}(t_{2}))-V_{2}(a_{2}(z_{2}))\right)^{-\frac{s_{2}}{p}}v_{1}^{1-p^{\prime}}(y_{1})\\
			&\qquad\times v_{2}^{1-p^{\prime}}(y_{2})dy_{1}dy_{2}\\
			&\quad+\left.\left.\prod_{i=1}^{2}\int_{b_{i}(t_{i})}^{b_{i}(z_{i})}\left(V_{i}(y_{i}))-V_{i}(a_{i}(z_{i}))\right)^{-\frac{s_{i}}{p}}v_{i}^{1-p^{\prime}}(y_{i})dy_{i}\right)^{q}u(z_{1},z_{2})dz_{1}dz_{2}\right)^{\frac{1}{q}}\\
			&\quad=\frac{p}{p-s_{1}}\frac{p}{p-s_{2}}\left(\int_{t_{1}}^{x_{1}}\int_{t_{2}}^{x_{2}}u(z_{1},z_{2})\prod_{i=1}^{2}\left(V_{i}(b_{i}(z_{i}))-V_{i}(a_{i}(z_{i}))\right)^{\frac{q(p-s_{i})}{p}}dz_{1}dz_{2}\right)^{\frac{1}{q}}.
			\end{split}
		}
		\end{equation}
		Combining $\eqref{eq:11}$, $\eqref{eq:2.18}$ together with $\eqref{eq:2.19}$ yields that
		$$
		{
			\begin{split}
		&\left(\int_{t_{1}}^{x_{1}}\int_{t_{2}}^{x_{2}}\left(V_{1}(b_{1}(z_{1}))-V_{1}(a_{1}(z_{1}))\right)^{\frac{q(p-s_{1})}{p}}
\left(V_{2}(b_{2}(z_{2}))-V_{2}(a_{2}(z_{2}))\right)^{\frac{q(p-s_{2})}{p}}u(z_{1},z_{2})dz_{1}dz_{2}\right)^{\frac{1}{q}}\\
		&\qquad\times\prod_{i=1}^{2}\frac{\frac{p}{p-s_{i}}}{\left(\left(\frac{p}{p-s_{i}}\right)^{p}+\frac{1}{s_{i}-1}\right)^{\frac{1}{p}}}\left(V_{i}(b_{i}(t_{i}))-V_{i}(a_{i}(x_{i}))\right)^{\frac{s_{i}-1}{p}}\leq C.
			\end{split}
			}
		$$
		Taking supermum over all $x_{i}$ and $t_{i}$ which satisfy $\eqref{eq:13}$ and applying $\eqref{eq:12}$, we have that the left hand side inequality in $\eqref{eq:14}$ holds.

		We now prove sufficiency. Assume that $\eqref{eq:12}$ holds. Here we will use some similar argument from \cite{HS,U}, and define two functions $a_{i}$ and $b_{i}$ for $i=1,2$
		be strictly increasing differentiable functions on $\mathbb{R}_{+}$ satisfying $\eqref{eq:2}$, then we obtain the inverse functions $a_{i}^{-1}$ and $b_{i}^{-1}$ are also exist, strictly increasing and differentiable. We define a sequence $\{m_{i}^{k}\}_{k\in\mathbb{Z}}$, $i=1,2$, as follow: for fixed $m_{i}>0$ write $m_{i}^{0}=m_{i}$ and
		\begin{equation}\label{eq:2.20}
		{
			\begin{split}
			&m_{i}^{k+1}=a_{i}^{-1}(b_{i}(m_{i}^{k})),\quad \text{if}\quad k\geq 0,\\
			&m_{i}^{k}=b_{i}^{-1}(a_{i}(m_{i}^{k+1})),\quad \text{if}\quad k\leq 0.
			\end{split}
			}
		\end{equation}
		Thus, we have
		\begin{equation}\label{eq:2.21}
		a_{i}(m_{i}^{k+1})=b_{i}(m_{i}^{k}),\quad \text{for all}\quad k\in\mathbb{Z}, i=1,2.
		\end{equation}
		
	In addition, we define the weight functions $u_{a_{1},a_{2}}$, $u_{b_{1},b_{2}}$, $u_{a_{1},b_{2}}$ and $u_{b_{1},a_{2}}$ by
		\begin{equation}\label{eq:2.22}
		{
			\begin{split}
			&u_{a_{1},a_{2}}(y_{1},y_{2})=u(a_{1}^{-1}(y_{1}),a_{2}^{-1}(y_{2}))(a_{1}^{-1})^{\prime}(y_{1})(a_{2}^{-1})^{\prime}(y_{2}),\\
			&u_{b_{1},b_{2}}(y_{1},y_{2})=u(b_{1}^{-1}(y_{1}),b_{2}^{-1}(y_{2}))(b_{1}^{-1})^{\prime}(y_{1})(b_{2}^{-1})^{\prime}(y_{2}),\\
			&u_{a_{1},b_{2}}(y_{1},y_{2})=u(a_{1}^{-1}(y_{1}),b_{2}^{-1}(y_{2}))(a_{1}^{-1})^{\prime}(y_{1})(b_{2}^{-1})^{\prime}(y_{2}),\\
			&u_{b_{1},a_{2}}(y_{1},y_{2})=u(b_{1}^{-1}(y_{1}),a_{2}^{-1}(y_{2}))(b_{1}^{-1})^{\prime}(y_{1})(a_{2}^{-1})^{\prime}(y_{2}),
			\end{split}
			}
		\end{equation}
		and
		$$a_{i}^{k}=a_{i}(m_{i}^{k}), b_{i}^{k}=b_{i}(m_{i}^{k}), \quad k\in\mathbb{Z}, i=1,2.$$
		Taking $g=fv_{1}^{1-p^{\prime}}v_{2}^{1-p^{\prime}}$ in the left hand side of $\eqref{eq:11}$, we obtain that it is less than or equal to
\begin{equation}\label{eq:2.23}
{
	\begin{split}
	&\left(\sum_{k_{1}\in\mathbb{Z}}\sum_{k_{2}\in\mathbb{Z}}\int_{a_{1}^{k_{1}}}^{b_{1}^{k_{1}}}\int_{a_{2}^{k_{2}}}^{b_{2}^{k_{2}}}
\left(\int_{y_{1}}^{b_{1}^{k_{1}}}\int_{y_{2}}^{b_{2}^{k_{2}}}g(t_{1},t_{2})dt_{1}dt_{2}\right)^{q}u_{a_{1},a_{2}}(y_{1},y_{2})dy_{1}dy_{2}\right)^{\frac{1}{q}}\\
	&\quad+\left(\sum_{k_{1}\in\mathbb{Z}}\sum_{k_{2}\in\mathbb{Z}}\int_{a_{1}^{k_{1}}}^{b_{1}^{k_{1}}}\int_{a_{2}^{k_{2}+1}}^{b_{2}^{k_{2}+1}}
\left(\int_{y_{1}}^{b_{1}^{k_{1}}}\int_{a_{2}^{k_{2}+1}}^{y_{2}}g(t_{1},t_{2})dt_{1}dt_{2}\right)^{q}u_{a_{1},b_{2}}(y_{1},y_{2})dy_{1}dy_{2}\right)^{\frac{1}{q}}\\
	&\quad+\left(\sum_{k_{1}\in\mathbb{Z}}\sum_{k_{2}\in\mathbb{Z}}\int_{a_{1}^{k_{1}+1}}^{b_{1}^{k_{1}+1}}\int_{a_{2}^{k_{2}}}^{b_{2}^{k_{2}}}
\left(\int_{a_{1}^{k_{1}+1}}^{y_{1}}\int_{y_{2}}^{b_{2}^{k_{2}}}g(t_{1},t_{2})dt_{1}dt_{2}\right)^{q}u_{b_{1},a_{2}}(y_{1},y_{2})dy_{1}dy_{2}\right)^{\frac{1}{q}}\\
	&\quad+\left(\sum_{k_{1}\in\mathbb{Z}}\sum_{k_{2}\in\mathbb{Z}}\int_{a_{1}^{k_{1}+1}}^{b_{1}^{k_{1}+1}}\int_{a_{2}^{k_{2}+1}}^{b_{2}^{k_{2}+1}}
\left(\int_{a_{1}^{k+1}}^{y_{1}}\int_{a_{2}^{k+1}}^{y_{2}}g(t_{1},t_{2})dt_{1}dt_{2}\right)^{q}u_{b_{1},b_{2}}(y_{1},y_{2})dy_{1}dy_{2}\right)^{\frac{1}{q}}\\
	&=II_{1}+II_{2}+II_{3}+II_{4}.
	\end{split}
	}		
\end{equation}	

Fix $t_{i}>0$ and let $t_{i}<z_{i}<x_{i}$, for $i=1,2$. Write $y_{i}=a_{i}(x_{i})$ in $\eqref{eq:12}$ and make the variable transformation $a_{i}(z_{i})=r_{i}$, $i=1,2$. Then we have
\begin{equation}\label{ss}
{
	\begin{split}
	&B(s_{1},s_{2})/\left(V_{1}(b_{1}(t_{1})-V_{1}(y_{1})\right)^{\frac{s_{1}-1}{p}}\left(V_{2}(b_{2}(t_{2})-V_{2}(y_{2})\right)^{\frac{s_{2}-1}{p}}\\
	&\geq\left(\int_{t_{1}}^{a_{1}^{-1}(y_{1})}\int_{t_{2}}^{a_{2}^{-1}(y_{2})}u(z_{1},z_{2})\prod_{i=1}^{2}\left(V_{i}(b_{i}(z_{i}))-V_{i}(a_{i}(z_{i}))\right)^{\frac{q(p-s_{i})}{p}}dz_{1}dz_{2}\right)^{\frac{1}{q}}\\
	&=\left(\int_{a_{1}(t_{1})}^{y_{1}}\int_{a_{2}(t_{2})}^{y_{2}}\prod_{i=1}^{2}\left(V_{i}(b_{i}(a_{i}^{-1}(r_{i})))-V_{i}(r_{i})\right)^{\frac{q(p-s_{i})}{p}}u_{a_{1},a_{2}}(r_{1},r_{2})dr_{1}dr_{2}\right)^{\frac{1}{q}}\\
	&\geq\left(\int_{a_{1}(t_{1})}^{y_{1}}\int_{a_{2}(t_{2})}^{y_{2}}\prod_{i=1}^{2}\left(V_{i}(b_{i}(t_{i})-V_{i}(r_{i})\right)^{\frac{q(p-s_{i})}{p}}u_{a_{1},a_{2}}(r_{1},r_{2})dr_{1}dr_{2}\right)^{\frac{1}{q}},
	\end{split}
	}
\end{equation}
where $u_{a_{1},a_{2}}$ is defined by $\eqref{eq:2.22}$. In the last estimate of $\eqref{ss}$ we  have used that
$$b_{1}(t_{1})<b_{1}(a_{1}^{-1}(r_{1})),\quad b_{2}(t_{2})<b_{2}(a_{2}^{-1}(r_{2}))$$

{\noindent}due to
$$t_{1}<z_{1}=a_{1}^{-1}(r_{1}), \quad t_{2}<z_{2}=a_{2}^{-1}(r_{2}).$$

Therefore, we conclude that
\begin{equation}\label{eq:2.24}
{
	\begin{split}
	&\left(\int_{c_{1}}^{y_{1}}\int_{c_{2}}^{y_{2}}\left(V_{1}(d_{1})-V_{1}(r_{1})\right)^{\frac{q(p-s_{1})}{p}}\left(V_{2}(d_{2})-V_{2}(r_{2})\right)^{\frac{q(p-s_{2})}{p}}u_{a_{1},a_{2}}(r_{1},r_{2})dr_{1}dr_{2}\right)^{\frac{1}{q}}\\
	&\qquad\qquad\qquad\qquad\quad\times\left(V_{1}(d_{1})-V_{1}(y_{1})\right)^{\frac{s_{1}-1}{p}}\left(V_{2}(d_{2})-V_{2}(y_{2})\right)^{\frac{s_{2}-1}{p}}<B(s_{1},s_{2})<\infty,
	\end{split}
	}
\end{equation}
for all $(c_{i},d_{i})=(a_{i}(t_{i}),b_{i}(t_{i}))$, $i=1,2$, in particular for all $(c_{i},d_{i})=(a_{i}^{k},b_{i}^{k})$. Hence, the condition $\eqref{eq:23}$ is satisfied. We can use Lemma $\ref{lem:3}$ again, and   we show that
\begin{equation}\label{eq:2.25}
{
	\begin{split}
	II_{1}&\leq\left(\sum_{k_{1}\in\mathbb{Z}}\sum_{k_{2}\in\mathbb{Z}}C^{q}\left(\int_{a_{1}^{k_{1}}}^{b_{1}^{k_{1}}}
\int_{a_{2}^{k_{2}}}^{b_{2}^{k_{2}}}g^{p}(y_{1},y_{2})v_{1}(y_{1})v_{2}(y_{2})dy_{1}dy_{2}\right)^{\frac{q}{p}}\right)^{\frac{1}{q}}\\
	&\leq C\left(\sum_{k_{1}\in\mathbb{Z}}\sum_{k_{2}\in\mathbb{Z}}\int_{a_{1}^{k_{1}}}^{b_{1}^{k_{1}}}\int_{a_{2}^{k_{2}}}^{b_{2}^{k_{2}}}f^{p}(y_{1},y_{2})v_{1}^{1-p^{\prime}}v_{2}^{1-p^{\prime}}(y_{2})dy_{1}dy_{2}\right)^{\frac{1}{p}}\\
	&=C\left(\int_{0}^{\infty}\int_{0}^{\infty}f^{p}(y_{1},y_{2})v_{1}^{1-p^{\prime}}v_{2}^{1-p^{\prime}}(y_{2})dy_{1}dy_{2}\right)^{\frac{1}{p}}.
	\end{split}
	}
\end{equation}
According to $\eqref{eq:24}$ and $\eqref{eq:2.24}$, it implies that
\begin{equation}\label{eq:2.26}
{
	\begin{split}
	C\leq\inf_{s_{1},s_{2}\in (1,p)}\tilde{A}_{W}(s_{1},s_{2})\prod_{i=1}^{2}\left(\frac{p-1}{p-s_{i}}\right)^{\frac{1}{p^{\prime}}}\leq\inf_{s_{1},s_{2}\in (1,p)}B(s_{1},s_{2})\prod_{i=1}^{2}\left(\frac{p-1}{p-s_{i}}\right)^{\frac{1}{p^{\prime}}}.
	\end{split}
	}
\end{equation}

Next, we will estimate   $II_{2}$. Let $t_{1}$ and $x_{2}$ be fixed, write $y_{1}=a_{1}(x_{1})$, $y_{2}=b_{2}(t_{2})$ in $\eqref{eq:12}$ and make the variable transformation $a_{1}(z_{1})=r_{1}$, $b_{2}(z_{2})=r_{2}$. Similar to estimate of $II_{1}$, we find that
\begin{equation}\label{eq:2.27}
{
	\begin{split}
	&\left(\int_{c_{1}}^{y_{1}}\int_{y_{2}}^{d_{2}}\left(V_{1}(d_{1})-V_{1}(r_{1})\right)^{\frac{q(p-s_{1})}{p}}\left(V_{2}(r_{2})-V_{2}(c_{2})\right)^{\frac{q(p-s_{2})}{p}}u_{a_{1},b_{2}}(r_{1},r_{2})dr_{1}dr_{2}\right)^{\frac{1}{q}}\\
	&\qquad\qquad\qquad\qquad\quad\times\left(V_{1}(d_{1})-V_{1}(y_{1})\right)^{\frac{s_{1}-1}{p}}\left(V_{2}(y_{2})-V_{2}(c_{2})\right)^{\frac{s_{2}-1}{p}}<B(s_{1},s_{2})<\infty,
	\end{split}
}
\end{equation}
for all $(c_{i},d_{i})=(a_{i}(t_{i}),b_{i}(t_{i}))$, $i=1,2$, where $u_{a_{1},b_{2}}$ is defined by $\eqref{eq:2.22}$. Then, $\eqref{eq:26}$ holds. We can conclude from  Lemma $\ref{lem:4}$ and the inequality $\eqref{eq:26}$ that
\begin{equation}\label{eq:2.28}
{
	\begin{split}
	II_{2}&\leq\left(\sum_{k_{1}\in\mathbb{Z}}\sum_{k_{2}\in\mathbb{Z}}C^{q}\left(\int_{a_{1}^{k_{1}}}^{b_{1}^{k_{1}}}
\int_{a_{2}^{k_{2}+1}}^{b_{2}^{k_{2}+1}}g^{p}(y_{1},y_{2})v_{1}(y_{1})v_{2}(y_{2})dy_{1}dy_{2}\right)^{\frac{q}{p}}\right)^{\frac{1}{q}}\\
	&\leq C\left(\sum_{k_{1}\in\mathbb{Z}}\sum_{k_{2}\in\mathbb{Z}}\int_{a_{1}^{k_{1}}}^{b_{1}^{k_{1}}}\int_{a_{2}^{k_{2}+1}}^{b_{2}^{k_{2}+1}}f^{p}(y_{1},y_{2})v_{1}^{1-p^{\prime}}v_{2}^{1-p^{\prime}}(y_{2})dy_{1}dy_{2}\right)^{\frac{1}{p}}\\
	&=C\left(\int_{0}^{\infty}\int_{0}^{\infty}f^{p}(y_{1},y_{2})v_{1}^{1-p^{\prime}}v_{2}^{1-p^{\prime}}(y_{2})dy_{1}dy_{2}\right)^{\frac{1}{p}}.
	\end{split}
}
\end{equation}
Thus it follows from  $\eqref{eq:27}$ and $\eqref{eq:2.27}$ that
\begin{equation}\label{eq:2.29}
{
	\begin{split}
	C\leq\inf_{s_{1},s_{2}\in (1,p)}\tilde{A}_{W}^{*}(s_{1},s_{2})\prod_{i=1}^{2}\left(\frac{p-1}{p-s_{i}}\right)^{\frac{1}{p^{\prime}}}\leq\inf_{s_{1},s_{2}\in (1,p)}B(s_{1},s_{2})\prod_{i=1}^{2}\left(\frac{p-1}{p-s_{i}}\right)^{\frac{1}{p^{\prime}}}.
	\end{split}
}
\end{equation}
For $III_{3}$, we first let $x_{1}$, $t_{1}$ be fixed, write $y_{1}=b_{1}(t_{1})$, $y_{2}=a_{2}(x_{2})$ in $\eqref{eq:12}$ and make the variable transformation $b_{1}(z_{1})=r_{1}$, $a_{2}(z_{2})=r_{2}$. Similar to estimate for $II_{1}$ we obtain that
\begin{equation}\label{eq:2.30}
{
	\begin{split}
	&\left(\int_{y_{1}}^{d_{1}}\int_{c_{2}}^{y_{2}}\left(V_{1}(r_{1})-V_{1}(c_{1})\right)^{\frac{q(p-s_{1})}{p}}\left(V_{2}(d_{2})-V_{2}(r_{2})\right)^{\frac{q(p-s_{2})}{p}}u_{b_{1},a_{2}}(r_{1},r_{2})dr_{1}dr_{2}\right)^{\frac{1}{q}}\\
	&\qquad\qquad\qquad\qquad\quad\times\left(V_{1}(d_{1})-V_{1}(y_{1})\right)^{\frac{s_{1}-1}{p}}\left(V_{2}(y_{2})-V_{2}(c_{2})\right)^{\frac{s_{2}-1}{p}}<B(s_{1},s_{2})<\infty,
	\end{split}
}
\end{equation}
for all $(c_{i},d_{i})=(a_{i}(t_{i}),b_{i}(t_{i}))$, $i=1,2$, where $u_{a_{1},b_{2}}$ is defined by $\eqref{eq:2.22}$. Thus, $\eqref{eq:19}$ holds. We can deduce from Lemma $\ref{lem:2}$ and the inequality $\eqref{eq:19}$ that
\begin{equation}\label{eq:2.31}
{
	\begin{split}
	II_{3}&\leq\left(\sum_{k_{1}\in\mathbb{Z}}\sum_{k_{2}\in\mathbb{Z}}C^{q}
\left(\int_{a_{1}^{k_{1}+1}}^{b_{1}^{k_{1}+1}}\int_{a_{2}^{k_{2}}}^{b_{2}^{k_{2}}}g^{p}(y_{1},y_{2})v_{1}(y_{1})v_{2}(y_{2})dy_{1}dy_{2}\right)^{\frac{q}{p}}\right)^{\frac{1}{q}}\\
	&\leq C\left(\sum_{k_{1}\in\mathbb{Z}}\sum_{k_{2}\in\mathbb{Z}}\int_{a_{1}^{k_{1}+1}}^{b_{1}^{k_{1}+1}}\int_{a_{2}^{k_{2}}}^{b_{2}^{k_{2}}}f^{p}(y_{1},y_{2})v_{1}^{1-p^{\prime}}v_{2}^{1-p^{\prime}}(y_{2})dy_{1}dy_{2}\right)^{\frac{1}{p}}\\
	&=C\left(\int_{0}^{\infty}\int_{0}^{\infty}f^{p}(y_{1},y_{2})v_{1}^{1-p^{\prime}}v_{2}^{1-p^{\prime}}(y_{2})dy_{1}dy_{2}\right)^{\frac{1}{p}}.
	\end{split}
}
\end{equation}
Applying $\eqref{eq:20}$ and $\eqref{eq:2.30}$ we have
\begin{equation}\label{eq:2.32}
{
	\begin{split}
	C\leq\inf_{s_{1},s_{2}\in (1,p)}{A}_{W}^{*}(s_{1},s_{2})\prod_{i=1}^{2}\left(\frac{p-1}{p-s_{i}}\right)^{\frac{1}{p^{\prime}}}\leq\inf_{s_{1},s_{2}\in (1,p)}B(s_{1},s_{2})\prod_{i=1}^{2}\left(\frac{p-1}{p-s_{i}}\right)^{\frac{1}{p^{\prime}}}.
	\end{split}
}
\end{equation}
For estimate for $II_{4}$, we first let $x_{i}$ be fixed, write $y_{i}=b_{i}(t_{i})$ in $\eqref{eq:12}$, $i=1,2$ and make the variable transformation $b_{1}(z_{1})=r_{1}$, $b_{2}(z_{2})=r_{2}$. Similar to estimate for $II_{1}$ we show that
\begin{equation}\label{eq:2.33}
{
	\begin{split}
	&\left(\int_{y_{1}}^{d_{1}}\int_{y_{2}}^{d_{2}}\left(V_{1}(r_{1})-V_{1}(c_{1})\right)^{\frac{q(p-s_{1})}{p}}\left(V_{2}(r_{2})-V_{2}(c_{2})\right)^{\frac{q(p-s_{2})}{p}}u_{b_{1},b_{2}}(r_{1},r_{2})dr_{1}dr_{2}\right)^{\frac{1}{q}}\\
	&\qquad\qquad\qquad\qquad\quad\times\left(V_{1}(d_{1})-V_{1}(y_{1})\right)^{\frac{s_{1}-1}{p}}\left(V_{2}(y_{2})-V_{2}(c_{2})\right)^{\frac{s_{2}-1}{p}}<B(s_{1},s_{2})<\infty,
	\end{split}
}
\end{equation}
for all $(c_{i},d_{i})=(a_{i}(t_{i}),b_{i}(t_{i}))$, $i=1,2$ and where $u_{b_{1},b_{2}}$ is defined by $\eqref{eq:2.22}$. Thus, $\eqref{eq:16}$ holds. We can conclude from Lemma A and the inequality \eqref{eq:16} that
\begin{equation}\label{eq:2.34}
{
	\begin{split}
	II_{4}&\leq\left(\sum_{k_{1}\in\mathbb{Z}}\sum_{k_{2}\in\mathbb{Z}}C^{q}\left(\int_{a_{1}^{k_{1}+1}}^{b_{1}^{k_{1}+1}}
\int_{a_{2}^{k_{2}+1}}^{b_{2}^{k_{2}+1}}g^{p}(y_{1},y_{2})v_{1}(y_{1})v_{2}(y_{2})dy_{1}dy_{2}\right)^{\frac{q}{p}}\right)^{\frac{1}{q}}\\
	&\leq C\left(\sum_{k_{1}\in\mathbb{Z}}\sum_{k_{2}\in\mathbb{Z}}\int_{a_{1}^{k_{1}+1}}^{b_{1}^{k_{1}+1}}
\int_{a_{2}^{k_{2}+1}}^{b_{2}^{k_{2}+1}}f^{p}(y_{1},y_{2})v_{1}^{1-p^{\prime}}v_{2}^{1-p^{\prime}}(y_{2})dy_{1}dy_{2}\right)^{\frac{1}{p}}\\
	&=C\left(\int_{0}^{\infty}\int_{0}^{\infty}f^{p}(y_{1},y_{2})v_{1}^{1-p^{\prime}}v_{2}^{1-p^{\prime}}(y_{2})dy_{1}dy_{2}\right)^{\frac{1}{p}}.
	\end{split}
}
\end{equation}
Combining $\eqref{eq:17}$ with $\eqref{eq:2.33}$ yields that
\begin{equation}\label{eq:2.35}
{
	\begin{split}
	C\leq\inf_{s_{1},s_{2}\in (1,p)}{A}_{W}(s_{1},s_{2})\prod_{i=1}^{2}\left(\frac{p-1}{p-s_{i}}\right)^{\frac{1}{p^{\prime}}}\leq\inf_{s_{1},s_{2}\in (1,p)}B(s_{1},s_{2})\prod_{i=1}^{2}\left(\frac{p-1}{p-s_{i}}\right)^{\frac{1}{p^{\prime}}}.
	\end{split}
}
\end{equation}
Combining $\eqref{eq:2.23}$, $\eqref{eq:2.25}$, $\eqref{eq:2.28}$, $\eqref{eq:2.31}$ together with $\eqref{eq:2.34}$, we find that the left hand side of $\eqref{eq:11}$ can be dominated by
$$C_{p,q}\left(\int_{0}^{\infty}\int_{0}^{\infty}f^{p}(y_{1},y_{2})v_{1}(y_{1})v_{2}(y_{2})dy_{1}dy_{2}\right)^{\frac{1}{p}},$$
where the best constant $C_{p,q}$ in  $\eqref{eq:14}$ can be estimated by  $\eqref{eq:2.26}$, $\eqref{eq:2.29}$, $\eqref{eq:2.32}$ and $\eqref{eq:2.35}$. The proof is complete. $\hfill \Box $

{\centering\section{A scale of weight characterizations for the P\'{o}lya-Knopp operator with variable limits}}

Our main result in this Section reads:
\begin{thm}\label{main:2}
		Let $0<p\leq q<\infty$, $1<s_{i}<p$, $i=1,2$ and let $u$ and $v$ be weight function on $\mathbb{R}_{+}^{2}$. Then
		the inequality
		\begin{equation}\label{eq:3.1}
		\|(G_{2}f)\|_{L^{q}(\mathbb{R}_{+}^{2},u)}\leq C_{p,q}^{*}\|f\|_{L^{p}(\mathbb{R}_{+}^{2},v)}
		\end{equation}
		holds for some finite constant $C_{p,q}^{*}$ and for all measurable functions $f\geq 0$ if and only if
		\begin{equation}\label{eq:3.2}
		{
			\begin{split}
				D(s_{1},s_{2})&=\sup\left(\int_{t_{1}}^{x_{1}}\int_{t_{2}}^{x_{2}}w(y_{1},y_{2})\left(b_{1}(y_{1})-a_{1}(y_{1})\right)^{-\frac{qs_{1}}{p}}\left(b_{2}(y_{2})-a_{2}(y_{2})
\right)^{-\frac{qs_{2}}{p}}dy_{1}dy_{2}\right)^{\frac{1}{q}}\\
				&\qquad \times\left(b_{1}(t_{1})-a_{1}(x_{1})\right)^{\frac{s_{1}-1}{p}}\left(b_{2}(t_{2})-a_{2}(x_{2})\right)^{\frac{s_{2}-1}{p}}<\infty,
			\end{split}
			}
		\end{equation}
		where the supremum is taken over all $x_{i}$ and $t_{i}$ such that $\eqref{eq:13}$ holds and  $w(x_{1},x_{2})$ is defined by $\eqref{eq_w}$.
		Furthermore, if $C_{p,q}^{*}$ is the best possible constant  in $\eqref{eq:3.1}$, then
			\begin{equation}\label{eq:3.3}
			\sup_{s_{1},s_{2}>1}D(s_{1},s_{2})\prod_{i=1}^{2}\left(\frac{e^{s_{i}}(s_{i}-1)}{e^{s_{i}}(s_{i}-1)+1}\right)^{\frac{1}{p}}\leq C_{p,q}^{*}\leq4\inf_{s_{1},s_{2}\in (1,p)}D(s_{1},s_{2})\prod_{i=1}^{2}\left(\frac{p-1}{p-s_{i}}\right)^{\frac{1}{p^{\prime}}}.
			\end{equation}	
\end{thm}
{\noindent}$Proof$. Assum that $\eqref{eq:3.2}$ holds. If we use the definition  $w(x_{1},x_{2})$ as in Theorem $\ref{main:2}$ and let
$$g(x_{1},x_{2})=f^{p}(x_{1},x_{2})v(x_{1},x_{2})$$ in $\eqref{eq:3.1}$. Then the  inequality $\eqref{eq:3.1}$ is equivalent to the following inequality
\begin{equation}\label{eq:3.4}
{
	\begin{split}
	&\left(\int_{0}^{\infty}\int_{0}^{\infty}(G_{2}g)^{\frac{q}{p}}(x_{1},x_{2}) w(x_{1},x_{2})dx_{1}dx_{2}\right)^{\frac{1}{q}}\leq C_{p,q}^{*}\left(\int_{0}^{\infty}\int_{0}^{\infty}g(x_{1},x_{2})dx_{1}dx_{2}\right)^{\frac{1}{p}}.
	\end{split}
	}
\end{equation}
Now replace $g(x_{1},x_{2})$ by $f^{p}(x_{1},x_{2})$, we have that $\eqref{eq:3.4}$ is equivalent to
\begin{equation}\label{eq:3.5}
{
	\begin{split}
	&\left(\int_{0}^{\infty}\int_{0}^{\infty}(G_{2}f)^{q}(x_{1},x_{2}) w(x_{1},x_{2})dx_{1}dx_{2}\right)^{\frac{1}{q}}\leq C_{p,q}^{*}\left(\int_{0}^{\infty}\int_{0}^{\infty}f^{p}(x_{1},x_{2})dx_{1}dx_{2}\right)^{\frac{1}{p}}.
	\end{split}
}
\end{equation}
Obviously  $\eqref{eq:3.1}$  is also equivalent to $\eqref{eq:3.5}$.

Moreover, it implies from Jensen's inequality  that
\begin{equation}\label{eq:3.6}
{
	\begin{split}
&\left(\int_{0}^{\infty}\int_{0}^{\infty}\left((G_{2}f)(x_{1},x_{2})\right)^{q} w(x_{1},x_{2})dx_{1}dx_{2}\right)^{\frac{1}{q}}\leq \left(\int_{0}^{\infty}\int_{0}^{\infty}(H_{2}f)^{q}(x_{1},x_{2})w(x_{1},x_{2})dx_{1}dx_{2}\right)^{\frac{1}{q}}.
	\end{split}
}
\end{equation}
Now we merely need to consider the inequality
\begin{equation}\label{eq:3.7}
\left(\int_{0}^{\infty}\int_{0}^{\infty}(H_{2}f)^{q}(x_{1},x_{2})w(x_{1},x_{2})dx_{1}dx_{2}\right)^{\frac{1}{q}}\leq {C}_{p,q}^{*}\left(\int_{0}^{\infty}\int_{0}^{\infty}f^{p}(x_{1},x_{2})dx_{1}dx_{2}\right)^{\frac{1}{p}}.
\end{equation}
Comparing with $\eqref{eq:11}$ we find that $\eqref{eq:3.7}$ is just $\eqref{eq:11}$ as $v_{1}(x_{1})=v_{2}(x_{2})=1$, $u(x_{1},x_{2})=w(x_{1},x_{2})(b_{1}(x_{1})-a_{1}(x_{1}))^{-q}(b_{2}(x_{2})-a_{2}(x_{2}))^{-q}$. Hence, according to condition $\eqref{eq:12}$ in Theorem $\ref{main:1}$, $\eqref{eq:3.7}$ holds if
$$
{
	\begin{split}
	B(s_{1},s_{2})&=\sup\left(\int_{t_{1}}^{x_{1}}\int_{t_{2}}^{x_{2}}w(y_{1},y_{2})\left(b_{1}(y_{1})-a_{1}(y_{1})\right)^{-\frac{qs_{1}}{p}}\left(b_{2}(y_{2})-a_{2}(y_{2})\right)^{-\frac{qs_{2}}{p}}dz_{1}dz_{2}\right)^{\frac{1}{q}}\\
	&\qquad \times\left(b_{1}(t_{1})-a_{1}(x_{1})\right)^{\frac{s_{1}-1}{p}}\left(b_{2}(t_{2})-a_{2}(x_{2})\right)^{\frac{s_{2}-1}{p}}<\infty,
	\end{split}
	}
$$
where the supremun is taken over all $x_{i}$ and $t_{i}$ which satisfy $\eqref{eq:13}$ for $i=1,2$. In fact, this condition is just  the assumption $\eqref{eq:3.2}$. Therefore, $\eqref{eq:3.7}$ holds and then  by $\eqref{eq:3.6}$, $\eqref{eq:3.1}$ holds. Furthermore, if ${C}_{p,q}^{*}$ is the best possible constant  in $\eqref{eq:11}$, then we have
$$
 {C}_{p,q}^{*}\leq4\inf_{s_{1},s_{2}\in (1,p)}B(s_{1},s_{2})\prod_{i=1}^{2}\left(\frac{p-1}{p-s_{i}}\right)^{\frac{1}{p^{\prime}}}.
$$
 Since $B(s_{1},s_{2})=D(s_{1},s_{2})$, we immediately have
\begin{equation}\label{eq:3.8}
{C}_{p,q}^{*}\leq4\inf_{s_{1},s_{2}\in (1,p)}D(s_{1},s_{2})\prod_{i=1}^{2}\left(\frac{p-1}{p-s_{i}}\right)^{\frac{1}{p^{\prime}}}.
\end{equation}

Now we assume that $\eqref{eq:3.1}$ and $\eqref{eq:3.4}$ holds with some finite constant $C_{p,q}^{*}$. Let $t_{i}$ and $x_{i}$ be fixed numbers, such that $0<t_{i}<x_{i}<\infty$ and $a_{i}(x_{i})<b_{i}(t_{i})$, $i=1,2$. For any $z_{i}\in (t_{i},x_{i})$ it follows that
$$a_{i}(z_{i})\leq a_{i}(x_{i})<b_{i}(t_{i})<b_{i}(z_{i})\leq b_{i}(x_{i}),$$ for $i=1,2$.
Next we now consider the following test function
$$
{
	\begin{split}
	g(y_{1},y_{2})=&\prod_{i=1}^{2}\left(b_{i}(t_{i})-a_{i}(z_{i})\right)^{-s_{i}}\chi_{(a_{i}(z_{i}),b_{i}(t_{i})}(y_{i})\\
	&+e^{-s_{2}}\left(b_{1}(t_{1})-a_{1}(z_{1})\right)^{-s_{1}}\left(y_{2}-a_{2}(z_{2})\right)^{-s_{2}}\chi_{(a_{1}(z_{1}),b_{1}(t_{1})}(y_{1})\chi_{(b_{2}(t_{2}),b_{2}(z_{2})}(y_{2})\\
	&+e^{-s_{1}}\left(y_{1}-a_{1}(z_{1})\right)^{-s_{1}}\left(b_{2}(t_{2})-a_{2}(z_{2})\right)^{-s_{2}}\chi_{(b_{1}(t_{1}),b_{1}(z_{1})}(y_{1})\chi_{(a_{2}(z_{2}),b_{2}(t_{2})}(y_{2})\\
	&+\prod_{i=1}^{2}e^{-s_{i}}\left(y_{i}-a_{i}(z_{i})\right)^{-s_{i}}\left(b_{i}(t_{i})-a_{i}(z_{i})\right)^{-s_{i}}\chi_{(b_{i}(t_{i}),b_{i}(z_{i})}(y_{i}).
	\end{split}
	}
$$
Similar to $\eqref{eq:2.18}$, the right hand side of $\eqref{eq:3.4}$ is less than or equal to
\begin{equation}\label{eq:3.9}
{
	\begin{split}
	\prod_{i=1}^{2}\left(1+e^{-s_{i}}\frac{1}{s_{i}-1}\right)(b_{i}(t_{i})-a_{i}(x_{i}))^{1-s_{i}}.
	\end{split}
	}
\end{equation}
For the left hand side in $\eqref{eq:3.4}$ we have
\begin{equation*}
{
	\begin{split}
	&II_{L}:=\left(\int_{0}^{\infty}\int_{0}^{\infty}(G_{2}g)^{\frac{q}{p}}(z_{1},z_{2}) w(z_{1},z_{2})dz_{1}dz_{2}\right)^{\frac{1}{q}}\geq\left(\int_{t_{1}}^{x_{1}}\int_{t_{2}}^{x_{2}}(G_{2}g)^{\frac{q}{p}}(z_{1},z_{2}) w(z_{1},z_{2})dz_{1}dz_{2}\right)^{\frac{1}{q}}\\
	&=\left(\int_{t_{1}}^{x_{1}}\int_{t_{2}}^{x_{2}}\left(\exp\left(\left(\prod_{i=1}^{2}\frac{1}{b_{i}(z_{i})-a_{i}(z_{i})}\right)\left[\int_{a_{1}(z_{1})}^{b_{1}(t_{1})}\int_{a_{2}(z_{2})}^{b_{2}(t_{2})}\ln
\left(\prod_{i=1}^{2}\left(b_{i}(t_{i})-a_{i}(z_{i})\right)^{-s_{i}}\right)dy_{1}dy_{2}\right.\right.\right.\right.\\
	&\qquad+\int_{a_{1}(z_{1})}^{b_{1}(t_{1})}\int_{b_{2}(t_{2})}^{b_{2}(z_{2})}\ln\left(e^{-s_{2}}(b_{1}(t_{1})-a_{1}(z_{1}))^{-s_{1}}(y_{2}-a_{2}(z_{2}))^{-s_{2}}\right)dy_{1}dy_{2}\\
	&\qquad+\int_{b_{1}(t_{1})}^{b_{1}(z_{1})}\int_{a_{2}(z_{2})}^{b_{2}(t_{2})}\ln\left(e^{-s_{1}}(y_{1}-a_{1}(z_{1}))^{-s_{1}}(b_{2}(t_{2})-a_{2}(z_{2}))^{-s_{2}}\right)dy_{1}dy_{2}\\
	&\qquad+\left.\left.\left.\left.\int_{b_{1}(t_{1})}^{b_{1}(z_{1})}\int_{b_{2}(t_{2})}^{b_{2}(z_{2})}\ln\left(\prod_{i=1}^{2}e^{-s_{i}}(y_{i}-a_{i}(z_{i}))^{-s_{i}}\right)dy_{1}dy_{2}\right]
\right)^{\frac{q}{p}}\right)w(z_{1},z_{2})dz_{1}dz_{2}\right)^{\frac{1}{q}}\\
	&=:\left(\int_{t_{1}}^{x_{1}}\int_{t_{2}}^{x_{2}}\left(\exp\left(\left(\prod_{i=1}^{2}\frac{1}{b_{i}(z_{i})-a_{i}(z_{i})}\right)[II_{1}+II_{2}+II_{3}+II_{4}]\right)^{\frac{q}{p}}\right)
w(z_{1},z_{2})dz_{1}dz_{2}\right)^{\frac{1}{q}},
	\end{split}
	}
\end{equation*}
	where
$$
{
	\begin{split}
	II_{1}&=-\left(\prod_{i=1}^{2}(b_{i}(t_{i})-a_{i}(z_{i}))\right)\sum_{i=1}^{2}s_{i}\ln(b_{i}(t_{i})-a_{i}(z_{i})),\\
	II_{2}&=-(s_{1}\ln(b_{1}(t_{1})-a_{1}(z_{1}))+s_{2})\int_{a_{1}(z_{1})}^{b_{1}(t_{1})}\int_{b_{2}(t_{2})}^{b_{2}(z_{2})}dy_{1}dy_{2}\\
&\quad-s_{2}\int_{a_{1}(z_{1})}^{b_{1}(t_{1})}\int_{b_{2}(t_{2})}^{b_{2}(z_{2})}\ln(y_{2}-a_{2}(z_{2}))dy_{1}dy_{2}\\
	&=-s_{2}(b_{1}(t_{1})-a_{1}(z_{1}))(b_{2}(z_{2})-a_{2}(z_{2}))\ln(b_{2}(z_{2})-a_{2}(z_{2}))\\
&\quad-s_{1}(b_{1}(t_{1})-a_{1}(z_{1}))(b_{2}(z_{2})-b_{2}(t_{2})\ln(b_{1}(t_{1})-a_{1}(z_{1}))\\
&\quad+s_{2}(b_{1}(t_{1})-a_{1}(z_{1}))(b_{2}(t_{2})-a_{2}(z_{2}))\ln(b_{2}(t_{2})-a_{2}(z_{2})),\\
	II_{3}&=-s_{1}(b_{1}(z_{1})-a_{1}(z_{1}))(b_{2}(t_{2})-a_{2}(z_{2}))\ln(b_{1}(z_{1})-a_{1}(z_{1}))\\
&\quad-s_{2}(b_{1}(z_{1})-b_{1}(t_{1}))(b_{2}(t_{2})-a_{2}(z_{2})) \ln(b_{2}(t_{2})-a_{2}(z_{2}))\\
&\quad+s_{1}(b_{1}(t_{1})-a_{1}(z_{1}))(b_{2}(t_{2})-a_{2}(z_{2}))\ln((b_{1}(t_{1})-a_{1}(z_{1}))\\
	\end{split}
	}
$$
and
$$
{
	\begin{split}
	II_{4}&=\sum_{i=1}^{2}-s_{i}\int_{b_{1}(t_{1})}^{b_{1}(z_{1})}\int_{b_{2}(t_{2})}^{b_{2}(z_{2})}dy_{1}dy_{2}-s_{i}\int_{b_{1}(t_{1})}^{b_{1}(z_{1})}\int_{b_{2}(t_{2})}^{b_{2}(z_{2})}\ln(y_{i}-a_{i}(z_{i}))dy_{1}dy_{2}\\
	&=-s_{1}(b_{1}(z_{1})-a_{1}(z_{1}))(b_{2}(z_{2})-b_{2}(t_{2}))\ln(b_{1}(z_{1})-a_{1}(z_{1}))\\
&\quad+s_{1}(b_{1}(t_{1})-a_{1}(z_{1}))(b_{2}(z_{2})-b_{2}(t_{2})) \ln(b_{1}(t_{1})-a_{1}(z_{1}))\\
&\quad-s_{2}(b_{1}(z_{1})-b_{1}(t_{1}))(b_{2}(z_{2})-z_{2}(z_{2}))\ln(b_{2}(z_{2})-z_{2}(z_{2}))\\
	&\quad+s_{2}(b_{1}(z_{1})-b_{1}(t_{1}))(b_{2}(t_{2})-a_{2}(z_{2}))\ln(b_{2}(t_{2})-a_{2}(z_{2})).
	\end{split}
	}
$$
Summing up we obtain
$$II_{1}+II_{2}+II_{3}+II_{4}=-\left(\prod_{i=1}^{2}(b_{i}(z_{i})-a_{i}(z_{i}))\right)\sum_{i=1}^{2}\ln(b_{i}(z_{i})-a_{i}(z_{i}))^{s_{i}}.$$
and we deduce that
\begin{equation}\label{eq:3.10}
{
	\begin{split}
	II_{L}&\geq\left(\int_{t_{1}}^{x_{1}}\int_{t_{2}}^{x_{2}}\left(\exp\left(-s_{1}\ln(b_{1}(z_{1})-a_{1}(z_{1}))-s_{2}\ln(b_{2}(z_{2})-a_{2}(z_{2}))\right)\right)^{\frac{q}{p}}w(z_{1},z_{2})dz_{1}dz_{2}\right)^{\frac{1}{q}}\\
	&=\left(\int_{t_{1}}^{x_{1}}\int_{t_{2}}^{x_{2}}(b_{1}(z_{1})-a_{1}(z_{1}))^{-\frac{s_{1}q}{p}}(b_{2}(z_{2})-a_{2}(z_{2}))^{-\frac{s_{2}q}{p}}w(z_{1},z_{2})dz_{1}dz_{2}\right)^{\frac{1}{q}}.
	\end{split}
}
\end{equation}
Thus, by combining $\eqref{eq:3.9}$ with $\eqref{eq:3.10}$ we obtain that
\begin{equation}\label{eq:3.11}
C_{p,q}^{*}\geq D(s_{1},s_{2})\prod_{i=1}^{2}\left(1+e^{-s_{i}}\frac{1}{s_{i}-1}\right)^{-\frac{1}{p}}=D(s_{1},s_{2})\prod_{i=1}^{2}\left(\frac{e^{s_{i}}(s_{i}-1)}{e^{s_{i}}(s_{i}-1)+1}\right)^{\frac{1}{p}}.
\end{equation}
Hence, $\eqref{eq:3.1}$ holds. Furthermore, by combining $\eqref{eq:3.8}$ with $\eqref{eq:3.11}$ we show that also $\eqref{eq:3.2}$ holds.

$\hfill \Box $

	\end{document}